\documentclass{amsart}
\usepackage[latin1]{inputenc}
\newtheorem{thm}{Theorem}

\newtheorem{prop}[thm]{Proposition}

\theoremstyle{definition}
\newtheorem{defn}[thm]{Definition}
\theoremstyle{remark}

\newcommand{\To}{\longrightarrow}

\begin{document}

\title{The number of weakly compact convex subsets of the Hilbert space}
\author{Antonio Avil\'es}
\address{University of Paris 7, Equipe de Logique
Math\'ematique, UFR de Mathématiques, 2 Place Jussieu, 75251
Paris, France} \email{aviles@logique.jussieu.fr,avileslo@um.es}

\begin{abstract}
We prove that for $\kappa$ an uncountable cardinal, there exist
$2^\kappa$ many non homeomorphic weakly compact convex subsets of
weight $\kappa$ in the Hilbert space $\ell_2(\kappa)$.
\end{abstract}

\thanks{The first author was supported by a Marie Curie
Intra-European Felloship MCEIF-CT2006-038768 and research projects
MTM2005-08379 (MEC and FEDER) and S\'{e}neca 00690/PI/04}


\subjclass[2000]{54B35; 52A07}

\maketitle

\section{Introduction}

The first example of two weakly compact convex subsets of a
Hilbert space which have the same uncountable weight but are not
homeomorphic is due to Corson and
Lindenstrauss~\cite{CorsonLindenstrauss}\cite{Lindenstrauss}, who
provided a nonseparable such set in which all points are
$G_\delta$, in contrast with a any closed ball. Examples of two
equivalent norms in a nonseparable Hilbert space whose closed
balls are not weakly homeomorphic was given in
\cite{AvilesHilbertBall}. Later in our joint work with O.
Kalenda~\cite{fiberorders} we developed the technique of fiber
orders, inspired on some Shchepin's ideas \cite{Shchepin1}, and we
applied it precisely to distinguish topologically a number of
weakly compact convex sets in the Hilbert space. We exhibited
there a countable family of nonhomeomorphic such spaces of a given
weight, and the natural question was posed to us by Gilles
Godefroy and Robert Deville about the cardinality of the set of
homeomorphism classes of these sets. In this note we show, using
again the machinery of fiber orders from~\cite{fiberorders} that
this cardinality is the greatest possible, namely $2^\kappa$ many
spaces of weight $\kappa$, for $\kappa$ an uncountable cardinal.
Let us remind that in the case of $\kappa=\omega$ it is a
consequence of Keller's theorem (cf. \cite{BessagaPelczynski} or
\cite{vanMill}) that there are $\omega$ many homeomorphism classes
of weakly compact convex subsets of $\ell_2(\omega)$, one for each
dimension, ranging from 0 to $\omega$.\\

The big family of nonhomeomorphic sets will be parametrized by a
family of trees. By a tree we mean a partially ordered set such
that every initial segment is well ordered and which has moreover
a minimum element, called the root of the tree. We shall write
$s\parallel t$ meaning that $s$ and $t$ are comparable elements of
a tree (i.e. either $s\leq t$ or $t\leq s$).\\

Given a tree $T$ and an uncountable set $\Gamma$ we consider the
compact space $K[T,\Gamma]\subset 2^{T\times\Gamma}$ to be the set
of all subsets $x\subset T\times \Gamma$ of cardinality at most 2
such that the first coordinates of elements of $x$ are contained
in a branch of $T$. The space $K[T,\Gamma]$ can be canonically
viewed as weakly compact subset $\tilde{K}$ of the Hilbert space
$\ell_2(K[T\times\Gamma])$, by identifying the empty set with $0$,
every singleton $\{(t,\gamma)\}\in K[T,\Gamma]$ with the
corresponding vector of the canonical euclidean base
$e_{\{(t,\gamma)\}}$, and every doubleton
$\{(t,\gamma),(t',\gamma')\}\in K[T,\Gamma]$ with the vector
$e_{\{(t,\gamma),(t,\gamma')\}}+e_{\{(t,\gamma)\}}+e_{\{(t,\gamma')\}}$.
Also, the convex hull of $\tilde{K}$ is then canonically
homeomorphic to the space of probability measures
$P(K[T,\Gamma])$. Our plan is to exhibit, for every cardinal
$\kappa>\omega$ a family of $2^\kappa$ many different trees
$\{\Upsilon_\Lambda : \Lambda<2^\kappa\}$ for which the
corresponding spaces $P(K[\Upsilon_\Lambda,\kappa])$ are
pairwise nonhomeomorphic.\\

\section{Preliminaries}

In this section, we briefly describe the technique of fiber orders developed in \cite{fiberorders}. We refer to this paper \cite{fiberorders} for more detailed explanations.\\

 Given a continuous surjection $f:K\To L$ and $y\in L$, a preorder relation is defined on the fiber $f^{-1}(y)$ by setting $x\leq x'$ if for every neighborhood $U$ of $x$ there exists a neighborhood $U'$ of $x'$ such that $f(U')\subset f(U)$. We denote by $\mathbb{F}_y(f)$ the set $f^{-1}(y)$ endowed with this preorder relation and by $\mathbb{O}_{y}(f)$ the associated ordered set obtained from $\mathbb{F}_y(f)$ modulo the equivalence relation given by $x\sim x'$ whenever $x\leq x'$ and $x'\leq x$.\\

Given $\tau$ an uncountable regular cardinal and $K$ a compact
space with $weight(K)\geq\tau$, we consider $\mathcal{Q}_\tau(K)$
the set of all quotients of $K$ of weight less than $\tau$,
endowed with its natural order: $(p_1:K\To L_1)\leq (p_2:K\To
L_2)$ if and only if there is a continuous surjection
$p_1^2:L_1\To L_2$ such that $p_1^2 p_1 = p_2$. In this case, we
call to the surjection $p_1^2:L_1\To L_2$ the natural internal
surjection between the two quotients. A subset
$\mathcal{S}\subset\mathcal{Q}_\tau(K)$ is called a
$\tau$-semilattice if the supremum of any family of less than
$\tau$ many elements of $\mathcal{S}$ is an element of
$\mathcal{S}$, and $\mathcal{S}$ is called cofinal if for every
$L\in\mathcal{Q}_\tau(K)$ there exists $L'\in \mathcal{S}$ such
that $L\leq L'$. A useful criterion is that a $\tau$-semilattice
$\mathcal{S}$ is cofinal if and only if every two different points
of $K$ are separated by some quotient in $\mathcal{S}$.\\

 Given a property $\mathcal{P}$, we say that the $\tau$-typical surjection of $K$ has property $\mathcal{P}$ if there is a cofinal $\tau$-semilattice $\mathcal{S}\subset\mathcal{Q}_\tau$ such that all natural internal surjection between elements of $\mathcal{S}$ have property $\mathcal{P}$, or equivalently if every cofinal $\tau$-semilattice $\mathcal{S}\subset \mathcal{Q}_\tau(K)$ contains a further cofinal $\tau$-semilattice such that all natural internal surjections between elements of $\mathcal{S}'$ have property $\mathcal{P}$.\\

In order to distinguish topologically two compact sets, we shall check that, for some uncountable regular cardinal $\tau$, the $\tau$-typical surjection of each one has fiber orders with different properties.\\

\section{Analysis of the fibers}

Given $S\subset T$ trees, and
$M\supset N$ infinite sets we consider the continuous surjection
$g:K[T\times M]\To K[S\times N]$ given by $g(x)= x\cap (S\times
N)$ and its associated surjection between spaces of probability
measures
$$f=P(g):P(K[T\times M])\To P(K[S\times N]).$$ The $\tau$-typical surjections in spaces
$P(K[\Upsilon,\Gamma])$ will have this form, so in this section we
shall analyze how the fiber orders of such a function $f$ look
like.\\

\subsection{The fibers of $g$}. According to the method
established in \cite{fiberorders} for the computation of fiber
orders in spaces of probability measures, the first step towards
the computation of the fiber orders of $f=P(g)$ is the analysis of
the fibers of $g$. For $y\in K[S\times N]$, there are three
cases:\\

Case 1: If $|y|=2$, then the fiber is trivial, $|g^{-1}(y)|=1$.\\

Case 2: If $|y|=1$, then there are two types of points in the
fiber of $y$. On the one hand we have the same $y\in g^{-1}(y)$,
around which $g$ is locally open, so $y$ is the minimum element of
$\mathbb{F}_y(g)$. On the other hand, if we take any other point
$x\in g^{-1}(y)$, then $|x|=2$ and $x$ is an isolated point, so it
has a neighborhood whose image is $\{y\}$, hence all those
elements are equivalent and maximum elements of the fiber. Hence, $\mathbb{O}_y(g)\cong\{0,1\}$.\\

Case 3: Finally, for $y=0$ we get the most interesting fiber. For
$x\in g^{-1}(0)$ we have the following three possibilities:\\

\begin{itemize}
\item If $x=0$, then $g$ is locally open around $x$, so $x$ is the
minimum element of the fiber.\\

\item If $|x|=2$, then $x$ is an
isolated point, so $x$ is a maximum element of the fiber.\\

\item If $|x|=1$, then $x=\{(t,\gamma)\} = [t,\gamma]$ is a
singleton with $(t,\gamma)\in T\times M^\setminus (S\times N)$. It
happens then that the image of a basic neighborhood of $x$ is a
set of the form
$$\{0\}\cup\{[s,\delta] : (s,\delta)\in (S\times N)\setminus F, s\parallel t\}$$ where $F$ is a finite
set.\\

Thus $[t,\gamma] \leq [t',\gamma']$ iff $\{s\in S : s\parallel
t\}\supset \{s\in S : s\parallel t'\}$.\\

In order to understand the structure of $\mathbb{O}_0(g)$, let us
call $R = R(S,T)$ the family of all subsets of $S$ of the form
$\{s\in S : s\parallel t\}$ for some $t\in T$, ordered by reverse
inclusion. This set is order isomorphic to the singletons of the
fiber of $0$, and the whole $\mathbb{O}_0(g)$ is obtained by
adding one mininum and one
maximum $\mathbb{O}_0(g) \cong \{-\infty\}\cup R\cup\{+\infty\}$.\\

\end{itemize}

\subsection{The fibers of $P(g)$}. It follows from
\cite{fiberorders} that for a surjection between spaces of
probability measures on scattered spaces of the form
$f=P(g):P(K)\To P(L)$, and $\mu = \sum_{i\in
I}\lambda_i\delta_{y_i}\in P(L)$, $\lambda_i>0$, we have that
$$\mathbb{O}_\mu(P(g))\cong\prod_{i\in
I}\mathbb{O}_{\delta_{y_i}}(P(g)),$$

where the product of ordered sets is endowed with the order
$(x_i)_{i\in I}\leq (y_i)_{i\in I}$ if and only if $x_i\leq y_i$
for all $i\in I$. Moreover, the fiber orders corresponding to a
Dirac measure $\delta_{y}\in P(L)$ is computed using the structure
of the fiber $g^{-1}(y)$, namely for $\nu,\nu'\in
P(g)^{-1}(\delta_y)$, we have that $\nu\leq\nu'$ if and only if
$$\nu\langle x_1,\ldots,x_n\rangle \leq \nu'\langle
x_1,\ldots,x_n\rangle$$ for every $x_1,\ldots,x_n\in g^{-1}(y)$.
The operation $\langle\cdot\rangle$ associates a subset of the
fiber $g^{-1}(y)$ to every finite subset of $g^{-1}(y)$ and it is
related to the fiber order. We shall not need here the definition
nor the precise computation of this operation, it will be enough
to take into account the following two basic properties of the
$\langle\cdot\rangle$ operation:\\

\begin{itemize}
\item[(P.1)] For every $x\in g^{-1}(y)$, $\langle x \rangle = \{z : z\geq x\}.$\\

\item[(P.2)] For every $x_1,\ldots,x_n\in g^{-1}(y)$, $\langle
x_1,\ldots,x_n\rangle$ is an upwards closed set, that is, if $z\in
\langle x_1,\ldots,x_n\rangle$ and $z'\geq z$, then $z'\in \langle
x_1,\ldots,x_n\rangle$.\\
\end{itemize}

We concentrate now on our map $f=P(g):P(K[T,M])\To P(K[S,N])$, and
on the computation of the fiber orders $\mathbb{O}_{\delta_y}(f)$.
From all the information above, there are already two cases where
this computation is clear:\\

Case 1: For $|y|=2$, $|\mathbb{O}_{\delta_y}(f)|=1$.\\

Case 2: For $|y|=1$, $\mathbb{O}_{\delta_y}(f)\cong [0,1]$.\\

We devote the rest of this section to the third case, when $y=0$.
Remember that we obtained that $\mathbb{O}_0(g)\cong \hat{R} =
\{-\infty\}\cup R\cup \{+\infty\}$, where $R$ is the family of sets of the form $\{s\in S : t\parallel s\}$ for $t\in T$, ordered by reverse inclusion. According to the results
mentioned above, $\mathbb{O}_{\delta_0}(f)$ can be described as
the set of all functions $\nu:\hat{R}\To [0,1]$ with $\nu(\hat{R})
= 1$ and endowed with the order that $\nu\leq\nu'$ if and only if
$$\nu\langle x_1,\ldots,x_n\rangle \leq \nu'\langle
x_1,\ldots,x_n\rangle$$ for every $x_1,\ldots,x_n\in\hat{R}$
(here, for a set $A\subset \hat{R}$, $\nu(A)$ denotes $\sum_{a\in
A}\nu(a)$).\\

We shall not give a complete description of what this ordered set
$\mathbb{O}_{\delta_0}(f)$ is, but we shall get from this just
some information which is relevant for us, namely about what we call linear walks on $\mathbb{O}_{\delta_0}(f)$.\\

\begin{defn}
Let $O$ be an ordered set and $t\in O$.\begin{enumerate} \item An element $s\in O$ is called an immediate successor of $t$ if $t<s$ and there is no $r\in O$ such that $t<r<s$. The set of immediate successors of $t$ in $O$ is called $imsuc_O(t)$.
\item An element $s\in O$ is called an linear successor of $t$ if
$s$ is a
maximal element of the set of all $x>t$ such that $\{y :
t\leq y \leq x\}$ is linearly ordered. The set of linear successors of $t$ in $O$ is called $lisuc_O(t)$.\\
\end{enumerate}

\end{defn}

\begin{defn}
Let $O$ be an ordered set and $\alpha$ an ordinal. A discrete walk of
length $\alpha$ in $O$ is a sequence $\{a_\beta :
\beta<\alpha\}\subset O$ fulfilling the following
conditions:\begin{enumerate} \item $a_0 = \min(O)$.
\item $a_{\beta+1}$ is an immediate successor of $a_\beta$, whenever $\beta+1<\alpha$.
\item $a_\beta=\sup\{a_\gamma : \gamma<\beta\}$ for every limit
ordinal $\beta<\alpha$.\\
\end{enumerate}
\end{defn}

\begin{defn}
Let $O$ be an ordered set and $\alpha$ an ordinal. A linear walk of
length $\alpha$ in $O$ is a sequence $\{a_\beta :
\beta<\alpha\}\subset O$ fulfilling the following
conditions:\begin{enumerate} \item $a_0 = \min(O)$,
\item $a_{\beta+1}$ is a linear successor of $a_\beta$, whenever $\beta+1<\alpha$. \item $a_\beta=\sup\{a_\gamma : \gamma<\beta\}$ for every limit
ordinal $\beta<\alpha$.\\
\end{enumerate}
\end{defn}

%

For every
$t\in \hat{R}$, we consider
$a[t]\in\mathbb{O}_{\delta_0}(f)$ given by $a[t](t) = 1$ and $a[t](s)=0$ for $s\neq t$.\\

\begin{prop}\label{meridianwalks}
Let $X=\mathbb{O}_{\delta_0}(f)$, $\{a_\beta :
\beta<\alpha\}\subset X$. The following are equivalent:
\begin{enumerate}
\item $\{a_\beta : \beta<\alpha\}$ is a linear walk on $X$. \item There exists $\{t_\beta : \beta<\alpha\}$ a discrete walk on $\hat{R}$,
such that $a_\beta = a[t_\beta]$ for every $\beta<\alpha$.\\
\end{enumerate}
\end{prop}

Proof: The statement of the proposition is a consequence of the
two facts stated below. All along the proof, in order to show that
$x\leq y$ in $X$ we use property (P.2) of the operation
$\langle\cdot\rangle$, so we check that $x(A)\leq y(A)$ for every
upwards closed set $A\subset \hat{R}$. The other way around, to
show that $x\not\leq y$ we use property (P.1), checking that there
is $t\in\hat{R}$ such that $x\{s : s\geq t\} > y \{s : s\geq
t\}$.\\

\emph{Fact 1}: For $t\in \hat{R}$, the linear successors of $a[t]$
are exactly the elements of the form $a[s]$, for $s$ immediate
successor of $t$. \emph{Proof}: First notice that if $y\leq
x=a[t]$, then $y\{s:s\geq t\}\geq x\{s:s\geq t\}=1$, so $y$ is
concentrated on the nodes of $\hat{R}$ which are greater or equal
than $t$. We claim that for $y>x$, the set $\{z:x\leq z\leq y\}$
is linearly ordered if and only if $y$ is of the form
$$y=y_{s,\lambda} = \lambda a[s] + (1-\lambda)a[t]$$ for some
$\lambda>0$ and some immediate successor $s$ of $t$. Observe that
these elements are ordered as $y_{s,\lambda}<y_{s',\lambda'}$ if
and only if $s=s'$ and $\lambda<\lambda'$, and therefore Fact 1
follows from this claim. Assume that $y>x$ is not of the form
$y_{s,\lambda}$, then either there is $p>t$ not immediate
successor of $t$ such that $y(t)=\lambda>0$ or otherwise there are
two immediate successors $s\neq s'$ of $t$ with $y(s)=\mu>0$ and
$y(s')=\mu'>0$. In the first case, we can find $t<p'<p$ and then
\begin{eqnarray*}
\frac{\lambda}{2}a[p] &+& (1-\frac{\lambda}{2})a[t]\ \text{ and}\\
\lambda a[p'] &+& (1-\lambda)a[t]
\end{eqnarray*}
are incomparable elements below $y$. In the second case,
\begin{eqnarray*}
\mu a[s] &+& (1-\mu)a[t]\ \text{ and}\\
\mu' a[s'] &+& (1-\mu')a[t]
\end{eqnarray*}
are incomparable elements below $y$. It remains to show that in fact a set of the form $\{z : x< z\leq y_{s,\lambda}\}$ is linearly ordered, but one can easilycheck that this set equals $\{y_{s,\mu} : \mu\leq \lambda\}$.\\

\emph{Fact 2}: If $s\in\hat{R}$ and $I=\{t_\gamma : \gamma<\beta\}$ is an increasing sequence of elements of $\hat{R}$ for some limit ordinal $\beta$, then $s=\sup\{t_\gamma : \gamma<\beta\}$ if and only if $a[s]=\sup\{a[t_\gamma] : \gamma<\beta\}$. \emph{Proof}: If $x$ is an upper bound of the set $\{a[t_\gamma] : \gamma<\beta\}$, then $x$ is concentrated on the set of upper bounds of $I$ because $x\{t:t\geq t_\gamma\} \geq a[t_\gamma]\{t: t\geq t_\gamma\}=1$. This fact implies immediately that if $s=\sup(I)$, then $a[s]=\sup\{a[t_\gamma] : \gamma<\beta\}$. For the converse, suppose that $x$ was not concentrated on the supremum of $I$, then $\lambda=x(u)>0$ for some upper bound $u$ of $I$ which is not the supremum, which means that there is another upper bound $v$ with $u\not\leq v$, and in this case $a[v]$ is an upper bound of $\{a[t_\gamma] : \gamma<\beta\}$ which is not greater than $x$.$\qed$\\

An ordered set $O$ is called irreducible if it is not isomorphic
to any product of two ordered sets of cardinality greater than one. Recall that $r_0=\{s\in S : s\parallel 0\}$ represents the minimum element of $R$, and hence the ``second element'' of $\hat{R}$.\\

\begin{prop}
The ordered set $\mathbb{O}_{\delta_0}(f)$ is irreducible.
\end{prop}

Proof: Suppose $\mathbb{O}_{\delta_0}(f)= X\times Y$. For every
$\lambda\in [0,1]$ consider the element $u_\lambda = \lambda a[r_0] +(1-\lambda) a[-\infty]$.\\

Since $\{v: v\leq u_1\} = \{u_\lambda : \lambda\in [0,1]\}$ is
linearly ordered, only one of the two coordinates of $u_1$ in
$X\times Y$ can be different from the minimum, say $u_1=(0,y_1)$,
and hence $u_\lambda=(0,y_\lambda)$, $y_\lambda\in Y$.\\

Notice that for every $v\in \mathbb{O}_{\delta_0}(f)$ different
from the minimum, there exists $\lambda>0$ with $v\geq u_\lambda$
(take $\lambda = 1-v(-\infty)$). But if we consider an element of
the form $(x,0)$, there is no $u_\lambda =(0,y_\lambda)$ below
this element except for $\lambda=0$. It follows that
$|X|=1$.$\qed$\\

Remember the result from \cite{fiberorders} that we mentioned at the beginning, that $\mathbb{O}_\mu(P(g))\cong\prod_{i\in
I}\mathbb{O}_{\delta_{y_i}}(P(g))$, for $\mu = \sum_{i\in
I}\lambda_i\delta_{y_i}\in P(L)$, $\lambda_i>0$. It follows now that there are exactly three types of fiber orders of $f$ which are irreducible: isomorphic to a singleton, to the interval $[0,1]$ or to the special ordered set $\mathbb{O}_{\delta_0}(f)$.\\

\section{The family of compact convex sets}

For every ordinal $\alpha$ let $\Upsilon_\alpha$ be a tree with
the following properties:\\

\begin{itemize}
\item The root of $\Upsilon_\alpha$ is the ordinal $\alpha$.\\
\item $\Upsilon_\alpha$ is an ever branching tree, that is,
every element has at least two immediate succesors.\\
\item The height of $\Upsilon_\alpha$ equals $\omega\cdot \alpha$ (ordinal product, meaning concatenation of $\alpha$ many copies of the ordinal $\omega$).\\
 \item $\Upsilon_\alpha$ has a
branch
of length $\omega\cdot\alpha$.\\
\item For every $\beta<\alpha$, the $\beta$-th level of
$\Upsilon_\alpha$ has cardinality $|\beta|$. In particular
$|\Upsilon_\alpha|=\max(|\alpha|,\omega)$.\\
\item $\Upsilon_\alpha$ is a Hausdorff tree, that is, for every node $t$ which is not an immediate successor, we have that $t=\sup\{s : s<t\}$.\\
\end{itemize}

Let now $\kappa$ be an uncountable cardinal. To every subset
$A\subset \kappa$ of cardinality $\kappa$ we associate a tree
$\Upsilon^A$ with the following properties:\\

\begin{itemize}
\item The set of immediate succesors of the root 0 of $\Upsilon^A$
equals the set $A$.\\

\item For every $\alpha\in A$, the subtree $\{t\in\Upsilon^A :
t\geq \alpha\}$ equals the tree $\Upsilon_\alpha$.\\
\end{itemize}

\begin{thm}
If $A\neq B$, then $P(K[\Upsilon^A,\kappa])$ is not homeomorphic
to $P(K[\Upsilon^B,\kappa])$.
\end{thm}

Proof: Let us assume that there exists $\delta\in A\setminus B$.
Let $\tau$ be an uncountable regular cardinal such that
$|\delta|<\tau\leq\kappa$. We shall find a difference between
$P(K[\Upsilon^A,\kappa])$ and $P(K[\Upsilon^B,\kappa])$ by looking
at the $\tau$-typical surjection.\\

Let $C\in\{A,B\}$ and let $\mathcal{F}_C$ be the family of all
subtrees $T\subset\Upsilon^C$ such that:\\

\begin{itemize}

\item[(F.1)] $|T|<\tau$.\\

\item[(F.2)] If $\gamma\in C$ and $\gamma\leq\delta$, then $\Upsilon_\gamma\subset T$.\\

\item[(F.3)] If $\gamma\in C$, $\gamma>\delta$, and
$\Upsilon_\gamma\cap T\neq\emptyset$, then the whole first $\omega
\cdot (\delta+1)$ levels of the
tree $\Upsilon_\gamma$ are contained in $T$.\\

\end{itemize}

The family of quotients of
$P(K[\Upsilon^C,\kappa])$ of the form $P(K[T,M])$ with
$T\in\mathcal{F}_C$ and $M$ a subset of $\kappa$ of cardinality
less than $\tau$ constitutes a cofinal $\tau$-semilattice in
$\mathcal{Q}_\tau(P(K[\Upsilon^C,\kappa]))$. Here, $P(K[T,M])$ is
viewed as quotient of $P(K[\Upsilon^C,\kappa])$ via the surjection
$$f=P(g):P(K[\Upsilon^C,\kappa])\To P(K[T,M]),$$
where $g:K[\Upsilon^C,\kappa]\To K[T,M]$ is given by $g(x) = x\cap
(T\times M)$. The fact that this is a $\tau$-semilattice follows from the fact that $\mathcal{F}_C$ is closed under unions of less than $\tau$ many
elements, and the fact that is cofinal (that is, it separate points) follows from the fact that $\Upsilon^C = \bigcup\mathcal{F}_C$.\\

Hence the $\tau$-typical surjection of $P(K[\Upsilon^C,\kappa])$
is of the form studied in the previous sections, $f:P(K[T,M])\To
P(K[S,N])$, with moreover $T,S\in\mathcal{F}_C$. We look at the
only irreducible fiber order in this surjection which is not
linearly ordered : $\mathbb{O}_{\delta_0}(f)$.\\

We say that two (discrete or linear) walks on an ordered set strongly intersect if the
first three elements of the two walks are the same (notice that
the first two elements of a discrete walk in $\hat{R}$ are always $-\infty$ and $r_0$).
The following two statements
(a) and (b) establish the difference between
$P(K[\Upsilon^A,\kappa])$ and $P(K[\Upsilon^B,\kappa])$.\\

\begin{itemize}
\item [(a)] In the only irreducible non linearly ordered fiber
order of the $\tau$-typical surjection of
$P(K[\Upsilon^A,\kappa])$ there is a linear walk of length
$\omega\cdot\delta$ which does not strongly
intersect any walk of length $\omega\cdot(\delta+1)$.\\

\item [(b)] In the only irreducible non linearly ordered fiber
order of the $\tau$-typical surjection of
$P(K[\Upsilon^B,\kappa])$, any linear walk of length $\omega\cdot\delta$ strongly intersects some linear walk of length $\omega\cdot(\delta+1)$.\\
\end{itemize}

By Proposition~\ref{meridianwalks}, these statements about linear walks in $\mathbb{O}_{\delta_0}(f)$ are equivalent to the corresponding statements about discrete walks in $\hat{R}=R\cup\{-\infty,+\infty\}$. Recall that $R$ is defined as the family of
all subsets of $S$ of the form $r_t= \{s\in S : s\parallel t\}$
for $t\in T$, endowed with the order reverse to inclusion. We need some information about the ordered set $R$:\\

\begin{itemize}

\item[(R.1)] First, recall that $r_0 = S$ is the minimum of $R$.\\

\item[(R.2)] The set of immediate successors of $r_0$ in $R$ equals $\{r_\gamma : \gamma\in S\cap\kappa\}$. Namely, if $t\not\geq \gamma$ for any $\gamma\in S\cap\kappa$, then since $T\in\mathcal{F}_C$, it follows that $r_t = \{0\}$ is the maximum of $R$.

We call $R_\gamma = \{r\in R : r\geq r_\gamma\}$ for every $\gamma\in S\cap\kappa$.\\

\item[(R.3)] If $\gamma\in S\cap\kappa$ and $\gamma\leq\delta$, then $R_\gamma$ is order isomorphic either to $\Upsilon_\gamma$ or to the result of adding one maximum element to $\Upsilon_\gamma$. \emph{Proof}: Let $r_t$ be an element of $R$. If $\alpha\leq t$ for some $\alpha\in (S\cap \kappa)\setminus\{\gamma\}$, then $r_t\not\in R_\gamma$ because $\alpha\in r_t\setminus r_\gamma$. Hence, $$R_\gamma = \{r_t : t\geq\gamma\}\cup\{r_t : t\geq\alpha,\ \alpha\not\in S\cap\kappa\}.$$
The righthandside of this union may be empty or consist of one element $\{0\}$ which is the maximum of $R$. The lefthandside equals $\{r_t : t\in\Upsilon_\gamma\}$. Since $\gamma\leq\delta$, $\Upsilon_\gamma\subset S$ and since $\Upsilon_\gamma$ is ever branching, the map $t\mapsto r_t$ is one-to-one and it is indeed an order isomorphism.\\

\item [(R.4)] If $\gamma\in S\cap\kappa$, $\gamma>\delta$, then $R_\gamma$ contains a discrete walk of length $\omega\cdot\gamma$. \emph{Proof}: Since $S\in \mathcal{F}_C$, the whole $\omega\cdot(\gamma+1)$ levels of $\Upsilon_\gamma$ are contained in $S$, and we can find $\{s_\beta : \beta<\omega\cdot(\gamma+1)\}$ an initial segment of length $\omega\cdot(\gamma+1)$ of $\Upsilon_\gamma$ contained in $S$. Again, because the tree $\Upsilon_\gamma$ is ever branching the map $\beta\mapsto r_{s_\beta}$ is one-to-one. Moreover, $\{r_{s_\beta} : \beta<\omega\cdot(\gamma+1)\}$ is a discrete walk. If it was not a discrete walk, there should exist an element $r_t$ different from all $r_{s_\beta}$'s and such that $r_{s_\alpha}<r_t<r_{s_\beta}$ for some $\alpha<\beta$. But clearly, in our situation, if $r_t<r_{s_\beta}$, then $t=r_{s_\alpha}$ for some $\alpha<\beta$.\\

\end{itemize}

In the case when $C=A$, due to (R.3) a branch of length
$\omega\cdot\delta$ in $\Upsilon_\delta$ induces, supplementing
$r_0$ at the beginning, a discrete walk on $R$ of length
$\omega\cdot\delta$ which passes through $r_\delta$. This walk
cannot strongly intersect any walk of length
$\omega\cdot(\delta+1)$ because $\Upsilon_\delta$ does not contain
any segment of length greater than $\omega\cdot\delta$.\\

In the case when $C=B$, suppose we are given a walk of length $\omega\cdot\delta$. This walk cannot pass through any $r_\gamma$ with $\gamma<\delta$ due to (R.3), since $\Upsilon_\gamma$ has height $\omega\cdot\gamma$. Since $\delta\not\in B$, the walk must pass trough $r_\gamma$ for some $\gamma>\delta$. But then by (R.4) there is a walk of length $\omega\cdot(\delta+1)$ which also passes through $r_\gamma$, that is, it strongly intersects our given walk.\\


\begin{thebibliography}{1}

\bibitem{AvilesHilbertBall}
A.~Avil\'{e}s, \emph{The unit ball of the {H}ilbert space in its
weak topology}, Proc. Am. Math. Soc. \textbf{135} (2007),
833--836.

\bibitem{fiberorders} A. Avilés and O. F. K. Kalenda, \emph{Fiber orders and compact spaces of uncountable weight}, to appear in Fundamenta Mathematicae.

\bibitem{BessagaPelczynski} C. Bessaga, A. Pelczynski,
\emph{Selected topics in infinite dimensional topology}.
Monografie matematyczne. Tom 58. Warszawa: PWN - Polish Scientific
Publishers (1975).

\bibitem{CorsonLindenstrauss} H. H. Corson, J. Lindenstrauss, \emph{On weakly compact subsets of Banach spaces}, Proc. Amer. Math. Soc. \textbf{17} (1966), 407--412.

\bibitem{Lindenstrauss} J. Lindenstrauss, \emph{Weakly compact sets - Their topological properties and the Banach
spaces that they generate.} Symposium on Infinite-Dimensional
Topology (Louisiana State Univ., Baton Rouge, La., 1967), pp.
235--273. Ann. of Math. Studies, No. 69, Princeton Univ. Press,
Princeton, N. J., 1972.

\bibitem{Shchepin1} E. V. Shchepin,
\emph{Topology of limit spaces of uncountable inverse spectra}.
Russ. Math. Surv. \textbf{31} (1976), No.~5, 155--191; translation
from Uspekhi Mat. Nauk \textbf{31} (1976), No.~5, 191--226.

\bibitem{vanMill} J. van Mill, \emph{Infinite-dimensional topology. Prerequisites and introduction}.
North-Holland Mathematical Library, 43 (1989).

\end{thebibliography}
\end{document}